\documentclass[12pt]{article}

\usepackage{graphicx}  
\usepackage{amsmath}  
\usepackage{hyperref}
\usepackage{amssymb,amsthm}
\usepackage{color}
\usepackage{tikz}
\usepackage{pgfplots}

\newcommand{\R}{{\mathbb R}}\newcommand{\N}{{\mathbb N}}

\newcommand{\Z}{{\mathbb Z}}\newcommand{\C}{{\mathbb C}}

\definecolor{brilliantrose}{rgb}{1.0, 0.33, 0.64}

\let\epsilon\varepsilon
\let\theta\vartheta

\newtheorem{theorem}{Theorem}[section]\newtheorem{lemma}[theorem]{Lemma}

\newtheorem{remark}[theorem]{Remark}

\usepackage{filecontents}

\title{Validity of the stochastic Landau approximation 
for super-pattern forming systems with a spatial 1:3 resonance}
\author{Anna Logioti$^1$ \\
{\small
Institut f\"ur Analysis, Dynamik und Modellierung,} \\ {\small Universit\"at Stuttgart, Pfaffenwaldring 57, } \\ {\small 70569 Stuttgart, Germany}
\\ {\small email: anna.logioti@mathematik.uni-stuttgart.de}}

\begin{document}

\maketitle

\begin{abstract}
We consider a  Kuramoto-Shivashinsky  like equation close to the threshold of  instability
with 
additive white noise and spatially periodic boundary conditions
which simultaneously exhibit Turing bifurcations  with a spatial 1:3 resonance
of the critical wave numbers.
For the description of the bifurcating solutions
we derive a system of coupled stochastic Landau equations.
It is the goal of this paper to prove error estimates between the associated 
approximation obtained through this amplitude system and true solutions of the original system. 
The Kuramoto-Shivashinsky  like equation serves as a prototype model for so-called super-pattern forming
systems with quadratic nonlinearity and additive white noise.

\medskip

{\bf Keywords. } super-pattern forming
systems, stochastic Landau equations, error estimates, Kuramoto-Shivashinsky
\medskip

{\bf MSC Classification. } 35R60, 35Q56, 35B32

\end{abstract}


%

\section{Introduction}

We are interested in the validity of amplitude equations in case of so called super-pattern forming 
systems. The spatially homogeneous ground state of these systems becomes unstable via multiple  
Turing  instabilities which simultaneously occur  
at two strictly positive wave numbers $k_1$ and $k_2$.  
It is the purpose of this paper to handle the case of 
a spatial resonance $ k_2 = 3 k_1 $. 
In case of spatially periodic boundary conditions the bifurcation scenario of such pattern forming systems
has been analyzed for instance in \cite{PK00} by deriving a system of amplitude equations 
which coincides with the lowest order approximation of the reduced system on the associated 
center manifold. 
We refer to  \cite{PK00} for an
overview about the occurrence and relevance of this situation and the  related  physical literature. 

Here, we are interested in  the same  situation but with additional small additive white noise
in time and space.
We derive a system of coupled stochastic Landau equations for the description of the bifurcating solutions.
It is the goal of this paper to prove error estimates between the associated 
approximation obtained through this amplitude system and true solutions of the super-pattern forming
systems.
We refrain from greatest generality and restrict our analysis to a doubly unstable  Kuramoto-Shivashinsky (duKS) equation 
\begin{equation}
\partial_t u  =  -(k_1^2+\partial_x^2)^2 (k_2^2+\partial_x^2)^2 u + \alpha u + \partial_x (u^2) + \xi \label{org1}
\end{equation}
where $x, \, u=u(x,t) \in \mathbb{R}, \, t \geq 0 $, $ \alpha \in  \mathbb{R} $ denotes the bifurcation parameter and $ \xi= \xi(x,t) $ stands for the 
noise and will be specified below. 
Other possible models can be found in Section \ref{secdisc}.
In the deterministic case, i.e.  $ \xi = 0 $,
 the linearization around the spatially homogeneous steady state $u=0$ of  \eqref{org1} is solved by $ u(x,t) = e^{\mathrm{i}kx+\lambda(k) t} $ with
\begin{equation}
\begin{aligned}
\lambda_j(k)=-(k_1^2-k^2)^2 (k_2^2-k^2)^2+\alpha.
\end{aligned}
\end{equation}
We see that the trivial solution $ u=0$ is spectrally  stable for $ \alpha \leq  0 $ and becomes unstable for $  \alpha  >   0 $ at the wave numbers $ \pm  k_j $ for $ j=1,2 $.
Since we are interested in the dynamics close to this first instability we introduce the small bifurcation parameter $\alpha=\varepsilon^2$ with $0<\varepsilon\ll 1$.
Interesting dynamics occur in particular in case of a spatial 1:3 resonance, i.e. if $ k_2 = 3 k_1 $.
The duKS equation \eqref{org1} serves as a prototype model for such super-pattern forming
systems with quadratic nonlinearity and additive white noise.

We consider the duKS equation with $ 2 \pi / k_1 $-spatially periodic boundary conditions.
Validity results for the stochastic Landau approximation of pattern forming systems 
are rare. In \cite{BMPS01} such a result has been established 
for the stochastic Landau approximation of the Swift-Hohenberg equation with noise,
\begin{equation}\label{swiho}
\partial_t u  =  -(1+\partial_x^2)^2 u + \alpha u - u^3 + \xi.
\end{equation}
In case of $ 2 \pi  $-spatially periodic boundary conditions
with the ansatz
\begin{equation}\label{Lapp}
u(x,t) = \varepsilon A(T)e^{ix} +c.c. + \mathcal{O}(\varepsilon^3),
\end{equation}
where $ \alpha = \varepsilon^2 $, $ T = \varepsilon^2 t \geq 0 $. and $ A(T) \in \C $,
a Landau equation
\begin{equation}
\partial_T A =  A - 3 A|A|^2 + \widetilde{\xi},
\end{equation}
with a rescaled noise 
$ \widetilde{\xi} = \widetilde{\xi}(X,T) $ has been derived and justified 
by proving error estimates for this approximation.
The situation when the pattern forming system contains not only cubic but also quadratic 
terms is more involved. 
An approximation result in this situation has been established for instance in 
\cite{B05}. The theory has been applied for instance to the 
Kuramoto-Shivashinsky (KS) equation 
\begin{equation} \label{KSeq}
\partial_{ t } u = - ( 1 + \partial_{ x }^{2})^{2} u + \alpha u + \partial_{x} (u^{2}) + \xi  ,
\end{equation}
and B\'enard’s problem.
In this paper we reconsider the validity question in case of quadratic terms again,  with the goal 
to bring the proof closer to the deterministic situation. This allows us  to transfer 
the analysis of more complicated deterministic bifurcation scenarios, such as the
simultaneous appearance of two Turing instabilities,
to the stochastic case, too,  in this  and in  future papers. 
Our slightly modified approach will add some new aspects to this problem, in particular 
we find an  amplitude system which is different from the one derived in \cite{B05} due to a different scaling of the noise terms.

The mathematically more challenging situation occurs when the spatial 
periodicity is dropped and if the problem is considered on the 
unbounded real line. Then for the Swift-Hohenberg equation \eqref{swiho}  
with the ansatz
\begin{equation}\label{GLapp}
u(x,t) = \varepsilon A(X, T)e^{ix} +c.c. + \mathcal{O}(\varepsilon^2),
\end{equation}
where $ T = \varepsilon^2 t \geq 0 $, $ X = \varepsilon x \in \R $, and $ A(X,T) \in \C $,
a Ginzburg-Landau equation
\begin{equation}
\partial_T A = 4 \partial_X^2 A + A -3  A|A|^2 + \widetilde{\xi},
\end{equation}
with a rescaled noise 
$ \widetilde{\xi} = \widetilde{\xi}(X,T) $ has been derived.
The Ginzburg-Landau equation appears as universal amplitude equation 
which can be derived whenever the original system exhibits a Turing instability.

In the deterministic case systems without quadratic nonlinearities 
can be handled trivially, cf. \cite{KSM92}, whereas the handling of 
systems with quadratic terms 
requires a lot of additional work, cf. \cite{vH91,Schn94,MS95,Schn99}.
Therefore, the situation for stochastic pattern systems with a quadratic nonlinearity is also more involved. 
The last decades have seen various attempts to justify the stochastic version
of this approximation by proving error estimates between the stochastic
Ginzburg-Landau approximations
and true solutions of the original systems,
see for instance \cite{BMPS01,BHP05,B07,BW13,BBS19}. 
Serious difficulties occur in doing so since 
stochastic PDEs with noise which is white in time and space
have to be solved.
The
focus of these papers was on the analytic handling of the noise terms and so 
in most cases the authors concentrate on toy problems, in particular on 
pattern forming systems
with cubic nonlinearities, more precisely, on systems without
quadratic nonlinearities. 
As a consequence, only a small part of the mathematical theory 
about the validity of the  Ginzburg-Landau approximation
developed  in the deterministic case has been transferred to the stochastic case so far. 
Therefore, by using an approach which is closer to the deterministic situation  we  strongly
expect that it is more easy also to  transfer these validity proofs to the stochastic situation.
However, the validity of the Ginzburg-Landau approximation is beyond the scope of this paper.

The plan of the paper is as follows.
In Section \ref{secsetup} we present the system which we will consider in detail, 
in particular we separate the solution $ u $ into an irregular purely stochastic part and a 
more regular part.
In Section \ref{Zsect} we discuss the analytic properties of the irregular purely stochastic part.
In Section \ref{secresult}
we derive the amplitude equation 
and formulate our  approximation result in Theorem 
\ref{mainth}.
In Section \ref{secerror} we introduce the error functions and derive the equations for the error.
In Section \ref{secresidual} we provide some estimates for the reduced residual terms.
In Section \ref{secerrorestimates} we  estimate the error by using Gronwall's inequality. We conclude this paper with a discussion given in 
Section \ref{secdisc}.
\medskip

{\bf Notation.} 
Possibly different constants which can be chosen independently of the small perturbation parameter 
$ 0 < \varepsilon \ll 1 $ are denoted with the same symbol $ C $. 
\medskip 

{\bf Acknowledgement.} I would like to thank especially Guido Schneider from Universität Stuttgart for introducing me to the analysis of stochastic amplitude equations for super-pattern forming systems. This work is partially supported by the Deutsche Forschungsgemeinschaft DFG through the cluster of excellence 'SimTech' under EXC 2075-390740016.

\section{Setup} \label{secsetup}

For notational simplicity we consider
the duKS equation \eqref{org1} with  critical wave numbers 
$ k_1 = 1 $ and $ k_2 = 3 $ as a prototype model.
Moreover, we choose 
$ 2 \pi $-spatially periodic boundary conditions
and a $ 2 \pi $-spatially periodic noise
\begin{equation}
\xi(x,t) = \sum_{k \in \Z} \alpha_{k} \widehat{\xi}(k,t) e^{ikx}.
\end{equation}
Herein,
$$ 
\widehat{\xi}(k,t) = \overline{\widehat{\xi}(-k,t)}   =  \widehat{\xi}_{r}(k,t) +  
\widehat{\xi}_{i}(k,t) ,
$$ 
where  $ \widehat{\xi}_{r}(k,\cdot) $ and $ \widehat{\xi}_{i}(k,\cdot) $
are distributional derivatives w.r.t. time $ t \geq 0 $ of standard Wiener processes 
$ \widehat{W}_{r}(k,\cdot) $ and $ \widehat{W}_{i}(k,\cdot) $
for each $ k \in \N_0 $. The coefficients 
$ \alpha_{k} = \alpha(k,\varepsilon) =  \overline{\alpha(-k,\varepsilon)} \in \C $ 
are specified below.
We also write  
$$ 
u(x,t) = \sum_{k \in \Z} \widehat{u}(k,t) e^{ikx},
$$  
with $ \widehat{u}(k,t) = \overline{\widehat{u}(-k,t)} $ 
due to the $ 2 \pi $-periodic boundary conditions.
With $ \lambda(k)  = -(1-k^2)^2 (9-k^2)^2 $
the Fourier coefficients satisfy 
\begin{eqnarray*}
\partial_t \widehat{u}(k,t) & = & \lambda(k) \widehat{u}(k,t) + \varepsilon^2 \widehat{u}(k,t)
+ i k (\widehat{u}*\widehat{u} )(k,t)
+ \alpha_k \widehat{\xi}(k,t), 
\end{eqnarray*}
where 
$$ (\widehat{u}*\widehat{v})(k,t) = \sum_{k' \in \Z} \widehat{u}(k-k',t)\widehat{v}(k',t)  $$ 
is the discrete convolution of $ \widehat{u} $ and $ \widehat{v} $.

Like for other systems with additive white noise we can separate 
$ \widehat{u} $ 
into an irregular purely stochastic part and a 
more regular part. We set 
$$ 
\widehat{u}(k,t) = \widehat{v}(k,t) +  \widehat{Z}(k,t),
$$ 
where the new variables satisfy 
\begin{equation}
\partial_t \widehat{v}(k,t) =  \lambda(k) \widehat{v}(k,t) + \varepsilon^2 \widehat{v}(k,t)
+ \varepsilon^2 \widehat{Z}(k,t)
+ i k ((\widehat{v} +  \widehat{Z})*(\widehat{v} +  \widehat{Z}))(k,t)
\end{equation}
and 
\begin{equation} \label{eqZ}
\partial_t \widehat{Z}(k,t) =  \lambda(k) \widehat{Z}(k,t)  + \alpha_k \widehat{\xi}(k,t).
\end{equation}
This separation has the advantage that $ \widehat{Z}(k,t) $ is at least 
continuous and $  \widehat{v}(k,t) $ at least differentiable w.r.t. time.
In the following, for each $k\in \mathbb N_0$, we choose initial conditons $  \widehat{Z}(k,0) = 0 $

\section{The processes $ \widehat{Z}(k,t) $} 
\label{Zsect}

We start with the equation  \eqref{eqZ} for $ \widehat{Z}(k,t) $ 
which  can be solved by the variation of constant formula.
Assuming $  \widehat{Z}(k,0) = 0 $ we find 
$$ 
\widehat{Z}(k,t) =  \int_0^t  e^{\lambda(k) (t-\tau)}  \alpha_k \widehat{\xi}(k,t) d\tau
= 
\int_0^t  e^{\lambda(k) (t-\tau)} \alpha_k d\widehat{W}(k,\tau)  ,
$$ 
where the integral is a stochastic integral in the It\^{o}-sense, cf. \cite{Oeksendal}.
$ \widehat{Z}(k,\cdot) $ is a so called  Ornstein-Uhlenbeck
process. These are well studied stochastic processes for which 
we recall some important estimates from the existing literature.
We have to distinguish the case $ k \neq \pm 1,\pm 3 $ and 
$ k = \pm 1,\pm 3 $.

{\bf a)} If $ \lambda(k) < 0 $ we use the following lemma  for estimating $ \widehat{Z}(k,\cdot)  $. 
\begin{lemma} \label{lem31}For $ k \neq \pm 1,\pm 3 $ we have  that
$$
P(\sup_{\tau \in [0,t]}|\widehat{Z}(k,\tau)| \geq c_k) \leq c_k^{-2}  \int_0^t e^{2 \lambda(k) (t-\tau)} |\alpha_k|^2  d \tau
 \leq \frac{|\alpha_k|^2}{2 |\lambda(k)| c_k^2}  .
$$ 
\end{lemma}
\noindent
{\bf Proof.}  This is a direct consequence of \cite[(5.1.4)]{Arnold}
where $ \widehat{Z}(k,t) = \int_0^t G(k,\tau) dW(k,\tau) $ is estimated by 
$$
P(\sup_{\tau \in [0,t]}|\widehat{Z}(k,\tau)| \geq c) \leq c^{-2}  \int_0^t E|G(k,\tau)|^2 d\tau.
$$
The statement of Lemma \ref{lem31} follows by setting $ G(k,\tau) = e^{\lambda(k) (t-\tau)} \alpha_k $.
\qed
\medskip

{\bf b)} If $ \lambda(k) = 0 $ we use the following lemma  for estimating $
\widehat{Z}(k,\cdot)  $. 
\begin{lemma} \label{lem32}
For $ k = \pm 1,\pm 3 $ we have  that
$$
P(\sup_{\tau \in [0,t]}|\widehat{Z}(k,\tau)| \geq c_k) \leq c_k^{-2}  \int_0^t |\alpha_k|^2  d \tau
 \leq \frac{|\alpha_k|^2 t}{c_k^2}  .
$$
\end{lemma}
\noindent
{\bf Proof.} 
We have 
$ \partial_t \widehat{Z}(k,t) =  \alpha_k \widehat{\xi}(k,t) $ 
and so $ \widehat{Z}(k,t) =  \alpha_k \widehat{W}(k,t) $.
Applying again \cite[(5.1.4)]{Arnold} as in Lemma \ref{lem31} but now with 
$ G(k,\tau) =  \alpha_k $ gives the statement of Lemma \ref{lem32}.
\qed
\medskip

In the following we need to know the order of $ \widehat{Z}(k,\tau) $ 
with respect to $ 0 < \varepsilon \ll 1 $ for all 
$ t \in [0,T_0/\varepsilon^2] $, i.e., on the natural time scale of the 
Landau approximation.
We choose $ c_k = \mathcal{O}(\varepsilon^2) $ for 
$ k \neq \pm 1,\pm 3 $ which, using Lemma \ref{lem31}, 
implies $ \alpha_k = \mathcal{O}(\varepsilon^2)$ 
for $ k \neq \pm 1,\pm 3 $ since the eigenvalues $ \lambda(k) $ are  independent of $ 0 < \varepsilon \ll 1 $.

For $ k =\pm 1,\pm 3 $ Lemma \ref{lem32} implies
$$
P(\sup_{\tau \in [0,t]}|\widehat{Z}(k,\tau)| \geq c_k)  \leq 
\mathcal{O}\bigg(\frac{|\alpha_k|^2}{c_k^2} \varepsilon^{-2}\bigg)
$$
and so if
we  choose $ c_k = \mathcal{O}(\varepsilon) $ for 
$ k =\pm 1,\pm 3 $ this  implies $ \alpha_k = \mathcal{O}(\varepsilon^2)$ for $ k = \pm 1,\pm 3 $.

For the subsequent error estimates we need an estimate on 
$ (\widehat{Z}(k,\tau))_{k \in \Z} $ in a suitable function space.
The error will be estimated w.r.t. the $ \ell^2_r $-norm 
$$ 
\| R \|_{\ell^2_r} = \big (\sum_{k \in \Z} |R_k|^2 (1+k^2)^r \big )^{1/2}
$$ 
for $ r \geq 1 $, and so we subsequently need an estimate for 
\begin{equation} \label{convul}
 \sup_{t \in [0,T_0/\varepsilon^2]} \| (\widehat{Z}_k(t) )_{k\in \Z} \|_{\ell^2_r}.
\end{equation}
Before we do so, we remark that the space $ \ell^2_r $ is closed 
under the discrete convolution if $ r > 1/2 $. In detail,
for $ r > 1/2 $
 there 
is a $ C > 0 $ such that for all $ \widehat{u},\widehat{v} \in \ell^2_r $ 
we have 
$$ 
\| \widehat{u} *\widehat{v}  \|_{\ell^2_r} \leq C \| \widehat{u}   \|_{\ell^2_r}\|  \widehat{v}  \|_{\ell^2_r}.
$$ 
By means of Lemma \ref{lem31} and setting $I_u = \{-3,-1,1,3\}$, we obtain 
\begin{eqnarray*}
\lefteqn{P(\sup_{t \in [0,T_0/\varepsilon^2]} \sum_{k\in \Z\setminus I_u}  |\widehat{Z}_k(t) |^2 (1+k^2)^r \geq C_{Z}^2) }
\\ & \leq & P( \sum_{k\in \Z\setminus I_u} \sup_{t \in [0,T_0/\varepsilon^2]} |\widehat{Z}_k(t) |^2 (1+k^2)^r \geq C_{Z}^2 C_s^{-2} \sum_{k \in \Z\setminus I_u} (1+k^2)^{-s}) 
\\ & \leq & \sum_{k\in \Z\setminus I_u} P(  \sup_{t \in [0,T_0/\varepsilon^2]} |\widehat{Z}_k(t) |^2 (1+k^2)^r \geq C_{Z}^2 C_s^{-2} (1+k^2)^{-s}) 
\\ & \leq & \sum_{k\in \Z\setminus I_u} P(  \sup_{t \in [0,T_0/\varepsilon^2]} |\widehat{Z}_k(t) |^2 \geq C_{Z}^2 C_s^{-2} (1+k^2)^{-s-r}) 
\\ & \leq & \sum_{k\in \Z\setminus I_u} P(  \sup_{t \in [0,T_0/\varepsilon^2]} |\widehat{Z}_k(t) | \geq C_{Z} C_s^{-1} (1+k^2)^{-(s+r)/2}) 
\\ & \leq & \frac{1}{C_{Z}^2}\sum_{k\in \Z\setminus I_u} \frac{|\alpha_k|^2  C_s^{2} (1+k^2)^{(s+r)}}{2 |\lambda(k)|} ,
\end{eqnarray*}
where $ C_s^2 = \sum_{k\in \Z\setminus I_u} (1+k^2)^{-s} < \infty $ for $ s > 1/2 $ and  where we used that 
$$ 
P(x+y\geq a+b)\leq P(x \geq a) + P(y \geq b)
$$ 
since $ \{(x,y): x+y\geq a+b \} \subset \{ (x,y) : x \geq a \vee y \geq b \} $.
Therefore, for a given $ \delta > 0 $ there exists a $ C_{Z} > 0 $ such that 
\begin{equation} \label{Zesti}
 P(\sup_{t \in [0,T_0/\varepsilon^2]} \sum_{k\in \Z\setminus I_u}  |\widehat{Z}_k(t) |^2 (1+|k|^2)^r \geq C_{Z}^2) 
\leq \delta 
\end{equation} 
if $ \alpha_k = \varepsilon^2 $ and 
\begin{equation} \label{Zesti1a}
\sum_{k\in \Z\setminus I_u} \frac{C_s^{2} (1+k^2)^{(s+r)} }{2| \lambda(k)|} = \mathcal{O}(1)< \infty 
\end{equation} 
which is the case if $ 2(s+r) < 7 $
or equivalently
$0 \leq  s+ r < 7/2 $
since $  \lambda(k) \sim -k^8 $ for $|k| \to \infty $.
Hence, the subsequent error equations can be solved in the space $ \ell^2_r $ with 
$ 1/2 < r < 3 $.

\section{Derivation of the amplitude equation and the approximation result}
\label{secresult}

We define the residual or formal error 
\begin{eqnarray*}
\textrm{Res}(k,t) 
& = & - \partial_t \widehat{v}(k,t) +  \lambda(k) \widehat{v}(k,t) + \varepsilon^2 \widehat{v}(k,t)
+ \varepsilon^2 \widehat{Z}(k,t) \\ && 
+ i k ((\widehat{v} +  \widehat{Z})*(\widehat{v} +  \widehat{Z}))(k,t)
\end{eqnarray*}
which counts how much a function $ \widehat{v} $ fails to satisfy the equations.
However, the full residual is not very useful for the following analysis
since various combinations of $ \widehat{Z} $ are better  
directly handled in the equations for the error.
Therefore, we will define a reduced residual below.

With these preparations we are now going to derive the
stochastic Landau approximation. 
We make the ansatz
$$
\widehat{v}(1,t) = \varepsilon A_1(\varepsilon^2 t), 
\qquad 
\widehat{v}(3,t) = \varepsilon A_3(\varepsilon^2 t),
$$
$$
\widehat{v}(-1,t) = \varepsilon A_{-1}(\varepsilon^2 t), 
\qquad 
\widehat{v}(-3,t) = \varepsilon A_{-3}(\varepsilon^2 t),
$$
for the unstable modes. 
For the derivation of the Landau equations 
we have to extend this ansatz to a few stable modes, too, namely 
$$
\widehat{v}(2,t) = \varepsilon^2 A_2(\varepsilon^2 t), 
\qquad 
\widehat{v}(4,t) = \varepsilon^2 A_4(\varepsilon^2 t),
\qquad 
\widehat{v}(6,t) = \varepsilon^2 A_{6}(\varepsilon^2 t), 
$$
$$
\widehat{v}(-2,t) = \varepsilon^2 A_{-2}(\varepsilon^2 t), 
\qquad 
\widehat{v}(-4,t) = \varepsilon^2 A_{-4}(\varepsilon^2 t),
\qquad 
\widehat{v}(-6,t) = \varepsilon^2 A_{-6}(\varepsilon^2 t).
$$
Moreover, we set   $$ \widehat{v}(k,t) = 0 $$ for $ \vert k \vert \not \in \{ 1,2,3,4,6\} $.
We set  $ \widehat{Z}(k,t) =c_k \widehat{Z}_k(t) $, such that $ \widehat{Z}_k(t) = \mathcal{O}(1) $.
We choose 
\begin{equation}\label{futureWork}
	 c_{\pm 1} = \varepsilon, \;  c_{\pm 3} = \varepsilon \quad \text{and} \quad c_k = \varepsilon^2  \text{ for }  k \neq \pm 1,\pm 3
\end{equation}and find
\begin{eqnarray*}
\textrm{Res}(0,t) & = & \mathcal{O}(\varepsilon^4),
 \\ 
\textrm{Res}(1,t) & = & \varepsilon^3 (- \partial_T A_1 + A_1 +  \widehat{Z}_1  
+ 2i (  \widehat{Z}_2  + A_2)
 ( \widehat{Z}_{-1}  + A_{-1})\\ &&
+ 2i \widehat{Z}_0
(  \widehat{Z}_1  + A_1) + 2i (  \widehat{Z}_3  + A_3)
 ( \widehat{Z}_{-2}  + A_{-2}) \\ && + 2i (  \widehat{Z}_4  + A_4)
 ( \widehat{Z}_{-3}  + A_{-3})
)+\mathcal{O}(\varepsilon^4), 
\\ 
\textrm{Res}(2,t) & = &
\varepsilon^2( \lambda(2) A_2 + 2i(  \widehat{Z}_1  + A_1)^2 
+ 4i(  \widehat{Z}_3  + A_3)( \widehat{Z}_{-1}  + A_{-1}) )
\\ && + \mathcal{O}(\varepsilon^4),
\\ 
\textrm{Res}(3,t) & = &  
 \varepsilon^3 (- \partial_T A_3 + A_3 +  \widehat{Z}_3  
+ 6i (  \widehat{Z}_2  + A_2)
 ( \widehat{Z}_{1}  + A_{1}) + 6i  \widehat{Z}_0  
 ( \widehat{Z}_{3}  + A_{3})
\\ &&  + 6i (  \widehat{Z}_6  + A_6)
 ( \widehat{Z}_{-3}  + A_{-3}) + 6i (  \widehat{Z}_4  + A_4)
 ( \widehat{Z}_{-1}  + A_{-1})
 )
+\mathcal{O}(\varepsilon^4),
\\ 
\textrm{Res}(4,t) & = & \varepsilon^2(\lambda(4) A_4 + 8i (  \widehat{Z}_3  + A_3)
 ( \widehat{Z}_{1}  + A_{1}))+ \mathcal{O}(\varepsilon^3),
\\ 
\textrm{Res}(5,t) & = &  \mathcal{O}(\varepsilon^3), 
 \\ 
\textrm{Res}(6,t) & = & \varepsilon^2( \lambda(6) A_6 + 6 i (  \widehat{Z}_3  + A_3)^2)
+ \mathcal{O}(\varepsilon^3),
\\ 
\textrm{Res}(k,t) & = &  \mathcal{O}(\varepsilon^3),
\end {eqnarray*}
for $ |k| \geq 7 $, and $ \textrm{Res}(-k,t) = \overline{\textrm{Res}(k,t)} $.
%
In order to make the  residual smaller we set 
\begin{eqnarray} \label{eqa1}
 \partial_T A_1 & = &  A_1 +  \widehat{Z}_1 + 2i (  \widehat{Z}_2  + A_2)
 ( \widehat{Z}_{-1}  + A_{-1})
+ 2i  \widehat{Z}_0
(\widehat{Z}_1  + A_1)
\\ && \nonumber
+ 2i (  \widehat{Z}_3  + A_3)
 ( \widehat{Z}_{-2}  + A_{-2})  + 2i (  \widehat{Z}_4  + A_4)
 ( \widehat{Z}_{-3}  + A_{-3})
, \\ 
0 & = & \lambda(2) A_2 + 2i(  \widehat{Z}_1  + A_1)^2 + 4i(  \widehat{Z}_3  + A_3)( \widehat{Z}_{-1}  + A_{-1}) , \label{A2eq} \\
 \partial_T A_3& = &  
 A_3 +  \widehat{Z}_3  
+ 6i (  \widehat{Z}_2  + A_2)
 ( \widehat{Z}_{1}  + A_{1}) + 6i  \widehat{Z}_0  
 ( \widehat{Z}_{3}  + A_{3})
\\ &&  + 6i (  \widehat{Z}_6  + A_6)
 ( \widehat{Z}_{-3}  + A_{-3}) + 6i (  \widehat{Z}_4  + A_4)
 ( \widehat{Z}_{-1}  + A_{-1}), \nonumber
\\ 
0 & = & \lambda(4) A_4 + 8i (  \widehat{Z}_3  + A_3)
 ( \widehat{Z}_{1}  + A_{1}), \label{A4eq}
 \\ 
0& = & \lambda(6) A_6 + 6 i (  \widehat{Z}_3  + A_3)^2. \label{A6eq}
\end{eqnarray}
We eliminate $ A_2 $, $ A_4 $, and $ A_6 $  in the equations for $ A_1 $ and $ A_3 $.
We  find
\begin{eqnarray} \label{SAE1}
 \partial_T A_1 & = &  A_1 +  \widehat{Z}_1 + 2i (  \widehat{Z}_2  - \frac{1}{\lambda(2)} (2i(  \widehat{Z}_1  + A_1)^2 \\ && \nonumber + 4i(  \widehat{Z}_3  + A_3)( \widehat{Z}_{-1}  + A_{-1}) ))
 ( \widehat{Z}_{-1}  + A_{-1})
 + 2i  \widehat{Z}_0
(\widehat{Z}_1  + A_1)
\\ && \nonumber
+ 2i (  \widehat{Z}_3  + A_3)
 ( \widehat{Z}_{-2}  - \frac{1}{\lambda(2)} ( -2i (\widehat{Z}_{-1}  + A_{-1})^2 - 4i(  \widehat{Z}_{-3}  + A_{-3})( \widehat{Z}_{1}  + A_{1})) )  
 \\ && \nonumber 
 + 2i (  \widehat{Z}_4  - \frac{8 i}{\lambda(4)}  (  \widehat{Z}_3  + A_3)( \widehat{Z}_{1}  + A_{1}))
 ( \widehat{Z}_{-3}  + A_{-3})
, \\  \label{SAE1neu}
 \partial_T A_3& = &  
 A_3 +  \widehat{Z}_3  
+ 6i (  \widehat{Z}_2  - \frac{1}{\lambda(2)} (2i(  \widehat{Z}_1  + A_1)^2 \\ && \nonumber + 4i(  \widehat{Z}_3  + A_3)( \widehat{Z}_{-1}  + A_{-1}) ))
 ( \widehat{Z}_{1}  + A_{1})
+ 6i  \widehat{Z}_0
(\widehat{Z}_3  + A_3
)
\\ &&  \nonumber + 6i (  \widehat{Z}_6  - \frac{6 i}{\lambda(6)}(  \widehat{Z}_3  + A_3)^2)
 ( \widehat{Z}_{-3}  + A_{-3}) \\ && \nonumber + 6i (  \widehat{Z}_4  - \frac{8 i}{\lambda(4)}  (  \widehat{Z}_3  + A_3)( \widehat{Z}_{1}  + A_{1}))
 ( \widehat{Z}_{-1}  + A_{-1}).
\end{eqnarray}
It is the goal of this paper to prove the following approximation result.
\begin{theorem} \label{mainth}
For all $ \delta > 0 $, $ C_1 > 0 $  there exist $ \varepsilon_0 > 0 $, $ C_2 > 0 $  such that for all $ \varepsilon \in (0,\varepsilon_0) $ the following holds.
Let $ (A_1,A_3) \in C([0,T_0],\C^2) $ be a solution of the system of stochastic amplitude equations
\eqref{SAE1}-\eqref{SAE1neu}
with 
\begin{equation} \label{ass3}
\sup_{T \in [0,T_0]} (|A_1(T)| + |A_3(T)|) \leq C_1.
\end{equation}
Then there are solutions 
$ u \in C([0,T_0/\varepsilon^2],H^1(\R,\R)) $
of the duKS equation  \eqref{KSeq} with 
$$ \mathcal{P}(\sup_{t \in [0,T_0/\varepsilon^2]} \sup_{x \in \R}|u(x,t)-(\varepsilon A_1(\varepsilon^2 t) e^{i x}+ \varepsilon A_3(\varepsilon^2 t) e^{3 i x}+ c.c) |  
\leq C_2 \varepsilon^2)   >   1 - \delta.
$$
\end{theorem}
\begin{remark}{\rm 
It is easy to see  that Assumption \eqref{ass3} can be  satisfied in the sense that for all $ \delta > 0 $ there exists a $ C_1 > 0 $ such that  
\begin{equation} \label{rem42}
\mathcal{P}(\sup_{T \in [0,T_0]} (|A_1(T)| + |A_3(T)|) \leq C_1)
> 1 - \delta.
\end{equation}
}
\end{remark}

\section{The equations for the error}
\label{secerror}

The error made by  the approximation from  
Section \ref{secresult}  is denoted by 
$ \varepsilon^2 R_{\pm 1} $ and $ \varepsilon^2 R_{\pm 3} $ and by
$ \varepsilon^3 R_{k} $ for $ k \not \in \{  \pm 1, \pm 3 \} $, respectively.
In detail, 
we set 
\begin{eqnarray*}
\widehat{v}(\pm 1,t) & = & \varepsilon A_{\pm 1}(\varepsilon^2 t)  + \varepsilon^2 R_{\pm 1}(t) , \\ 
\widehat{v}(\pm 3,t) & = & \varepsilon A_{\pm 3}(\varepsilon^2 t)  + \varepsilon^2 R_{\pm 3}(t) , \\ 
\widehat{v}(\pm 2,t) & = & \varepsilon^2 A_{\pm 2}(\varepsilon^2 t)   + \varepsilon^3 R_{\pm 2}(t) , \\ 
\widehat{v}(\pm 4,t) & = & \varepsilon^2 A_{\pm 4}(\varepsilon^2 t)   + \varepsilon^3 R_{\pm 4}(t) , \\ 
\widehat{v}(\pm 6,t) & = & \varepsilon^2 A_{\pm 6}(\varepsilon^2 t)   + \varepsilon^3 R_{\pm 6}(t) , 
\end{eqnarray*}
and $ \widehat{v}(k,t) = \varepsilon^3 R_{k}(t) $ for $ k = 0,\pm 5 $ and $ |k| \geq 7 $. 
Moreover, we use the notation $ Y_{\pm 2} = A_{\pm 2} + \widehat{Z}_{\pm 2} $,
$ Y_{\pm 4} = A_{\pm 4} + \widehat{Z}_{\pm 4} $, $ Y_{\pm 6} = A_{\pm 6} + \widehat{Z}_{\pm 6} $
and $ Y_{k} = \widehat{Z}_k $ for $ k \neq \pm 1,\pm 2, \pm 3, \pm 4, \pm 6  $. We further recall that $ \widehat{Z}(k,t) =\varepsilon\widehat{Z}_k(t) $ for $k=\pm 1, \pm 3$ and 
$ \widehat{Z}(k,t) =\varepsilon^2 \widehat{Z}_k(t)$ otherwise.
We find 
\begin{eqnarray*}
\partial_t (\varepsilon A_{1}  + \varepsilon^2 R_{1} )& = & \lambda(1)  (\varepsilon A_{1}  + \varepsilon^2 R_{1} )+ \varepsilon^2 (\varepsilon A_{1}  + \varepsilon^2 R_{1} )
+ \varepsilon^3 \widehat{Z}_1
\\ && +  i  \sum_{1- k' = \pm 1,\pm 3 }(\varepsilon A_{1-k'}  + \varepsilon^2 R_{1-k'} + \varepsilon  \widehat{Z}_{1-k'})(\varepsilon^3 R_{k'} + \varepsilon^2 Y_{k'})
\\ && + i  \sum_{k' \in I_1 \cap \vert k^\prime \vert \in \{1,3\}}(\varepsilon^3 R_{1-k'} + \varepsilon^2 Y_{1-k'})(\varepsilon A_{k'}  + \varepsilon^2 R_{k'} + \varepsilon  \widehat{Z}_{k'})
\\ && + i  \sum_{k' \in I_1 \cap \vert k^\prime \vert \notin \{1,3\}}(\varepsilon^3 R_{1-k'} + \varepsilon^2 Y_{1-k'})(\varepsilon^3 R_{k'} + \varepsilon^2 Y_{k'}),
\end{eqnarray*}
\begin{eqnarray*}
\partial_t (\varepsilon A_{3}  + \varepsilon^2 R_{3} )& = & \lambda(3)  (\varepsilon A_{3}  + \varepsilon^2 R_{3} )+ \varepsilon^2 (\varepsilon A_{3}  + \varepsilon^2 R_{3} )
+ \varepsilon^3 \widehat{Z}_3
\\ && + 3 i  \sum_{3- k' = \pm 1,\pm 3 }(\varepsilon A_{3-k'}  + \varepsilon^2 R_{3-k'} + \varepsilon  \widehat{Z}_{3-k'})(\varepsilon^3 R_{k'} + \varepsilon^2 Y_{k'})
\\ && + 3 i  \sum_{k' \in I_3 \cap \vert k^\prime \vert \in \{1,3\}}(\varepsilon^3 R_{3-k'} + \varepsilon^2 Y_{3-k'})(\varepsilon A_{k'}  + \varepsilon^2 R_{k'} + \varepsilon  \widehat{Z}_{k'})
\\ && + 3 i  \sum_{k' \in I_3 \cap \vert k^\prime \vert \notin \{1,3\}}(\varepsilon^3 R_{3-k'} + \varepsilon^2 Y_{3-k'})(\varepsilon^3 R_{k'} + \varepsilon^2 Y_{k'}),
\end{eqnarray*}
%
and
\begin{eqnarray*}
\lefteqn{\partial_t \varepsilon^3 R_{k} 
+ \partial_t \varepsilon^2 (\delta_{k,2} A_2 + \delta_{k,4} A_4 +\delta_{k,6} A_6)
+ \partial_t \varepsilon^2 (\delta_{k,-2} A_{-2}+ \delta_{k,-4} A_{-4}+\delta_{k,-6} A_{-6})}
\\
& = & \lambda(k) \varepsilon^3 R_{k} + \lambda(k)\varepsilon^2 (\delta_{k,2} A_2 + \delta_{k,4} A_4 +\delta_{k,6} A_6) +\lambda(k) \varepsilon^2 (\delta_{k,-2} A_{-2}+ \delta_{k,-4} A_{-4}+\delta_{k,-6} A_{-6})
  \\&& +\varepsilon^2 \varepsilon^3 R_{k} +\varepsilon^4(\delta_{k,2} A_2 + \delta_{k,4} A_4 +\delta_{k,6} A_6) +\varepsilon^4(\delta_{k,-2} A_{-2}+ \delta_{k,-4} A_{-4}+\delta_{k,-6} A_{-6})
+ \varepsilon^4 \widehat{Z}_k
\\&&+ i k \sum_{k' \in I_k \cap \vert k^\prime \vert \in \{1,3\}}(\varepsilon^3 R_{k-k'} + \varepsilon^2 Y_{k-k'})(\varepsilon A_{k'}  + \varepsilon^2 R_{k'} + \varepsilon  \widehat{Z}_{k'})
\\ && + i k \sum_{\vert k- k' \vert = 1,3\; \cap \; \vert k^\prime \vert \notin \{1,3\}}(\varepsilon A_{k-k'}  + \varepsilon^2 R_{k-k'} + \varepsilon  \widehat{Z}_{k-k'})(\varepsilon^3 R_{k'} + \varepsilon^2 Y_{k'})
\\ && + i k \sum_{\vert k- k' \vert = 1,3\; \cap \;\vert k^\prime \vert \in \{1,3\}}(\varepsilon A_{k-k'}  + \varepsilon^2 R_{k-k'} + \varepsilon  \widehat{Z}_{k-k'})(\varepsilon A_{k'}  + \varepsilon^2 R_{k'} + \varepsilon  \widehat{Z}_{k'})
\\ && + i k \sum_{k' \in I_k \cap \vert k^\prime \vert \notin \{1,3\}}(\varepsilon^3 R_{k-k'} + \varepsilon^2 Y_{k-k'})(\varepsilon^3 R_{k'} + \varepsilon^2 Y_{k'}),
\end{eqnarray*}
with $ I_k = \Z\setminus\{|k'|,|k-k'|  \in \{1,3\}\} $. 
Reordering the terms in the equations for $ R_1 $  shows that
the error function $ R_1 $ satisfies
\begin{equation} \label{r1eq}
\partial_t R_1 = \lambda(1) R_1 + f_1(R) ,
\end{equation}
with 
\begin{eqnarray*}
\varepsilon^2 f_1(R) 
& = &  \varepsilon^4 R_{1} 
\\ && + i  \sum_{1- k' = \pm 1,\pm 3 }\big(\varepsilon^2 R_{1-k'}(\varepsilon^3 R_{k'} +\varepsilon^2 Y_{k'}) +\varepsilon Y_{1-k'} \varepsilon^3 R_{k'} \big)
 \\&&+  i  \sum_{k' \in I_1 \cap \vert k^\prime \vert \in \{1,3\} } \big(\varepsilon^3 R_{1-k'} (\varepsilon^2 R_{k'}+\varepsilon Y_{k'}) +\varepsilon^2 Y_{1-k'} \varepsilon^2 R_{k'} \big)\\ && + i  \sum_{k' \in I_1 \cap \vert k^\prime \vert \notin \{1,3\}}(\varepsilon^3 R_{1-k'} + \varepsilon^2 Y_{1-k'})(\varepsilon^3 R_{k'} + \varepsilon^2 Y_{k'}) +  \textrm{Res}_r(1,t),
\end{eqnarray*}
and the reduced residual  
\begin{eqnarray*}
\textrm{Res}_r(1,t) & = & \varepsilon^3 (- \partial_T A_1 + A_1 +  \widehat{Z}_1  
+ 2i (  \widehat{Z}_2  + A_2)
 ( \widehat{Z}_{-1}  + A_{-1})\\ &&
+ 2i \widehat{Z}_0
(  \widehat{Z}_1  + A_1) + 2i (  \widehat{Z}_3  + A_3)
 ( \widehat{Z}_{-2}  + A_{-2}) \\ && + 2i (  \widehat{Z}_4  + A_4)
 ( \widehat{Z}_{-3}  + A_{-3})
)= 0 , 
\end{eqnarray*}
since $ A_1 $ satisfies \eqref{SAE1}. We have 
$ \textrm{Res}_r(1,t) = \textrm{Res}(1,t) $ due to the fact that the 
nonlinearity only contains quadratic terms.
Doing the same 
in the equations for $ R_3 $  shows that
the error function $ R_3 $ satisfies
\begin{equation} \label{r3eq}
\partial_t R_3 = \lambda(3) R_1 + f_3(R) ,
\end{equation}
with 
\begin{eqnarray*}
\varepsilon^2 f_3(R) 
& = &  \varepsilon^4 R_{3} 
\\ && + 3 i  \sum_{3- k' = \pm 1,\pm 3 }\big(\varepsilon^2 R_{3-k'}(\varepsilon^3 R_{k'} +\varepsilon^2 Y_{k'}) +\varepsilon Y_{3-k'} \varepsilon^3 R_{k'} \big)
\\&&+ 3 i  \sum_{k' \in I_3 \cap \vert k^\prime \vert \in \{1,3\}}\big(\varepsilon^3 R_{3-k'} (\varepsilon^2 R_{k'}+\varepsilon Y_{k'}) +\varepsilon^2 Y_{3-k'} \varepsilon^2 R_{k'} \big)
\\ && + 3 i  \sum_{k' \in I_3 \cap \vert k^\prime \vert \notin \{1,3\}}(\varepsilon^3 R_{3-k'} + \varepsilon^2 Y_{3-k'})(\varepsilon^3 R_{k'} + \varepsilon^2 Y_{k'})+  \textrm{Res}_r(3,t),
\end{eqnarray*}
and the reduced residual   
\begin{eqnarray*}
\textrm{Res}_r(3,t) & = &  
 \varepsilon^3 (- \partial_T A_3 + A_3 +  \widehat{Z}_3  
+ 6i (  \widehat{Z}_2  + A_2)
 ( \widehat{Z}_{1}  + A_{1})
\\ &&  + 6i (  \widehat{Z}_6  + A_6)
 ( \widehat{Z}_{-3}  + A_{-3}) + 6i (  \widehat{Z}_4  + A_4)
 ( \widehat{Z}_{-1}  + A_{-1})
 ) = 0
\end{eqnarray*}
since $ A_3 $ satisfies \eqref{SAE1neu}.
We have 
$ \textrm{Res}_r(3,t) = \textrm{Res}(3,t) $ due to the fact that the 
nonlinearity only contains quadratic terms.
Reordering the terms in the equations for $ R_k $ shows that
the error function $ R_k $  for $ |k|  \neq \pm 1,\pm 3 $ satisfy
\begin{equation} \label{rkeq}
\partial_t R_k = \lambda(k) R_k + f_k(R) 
\end{equation}
with 
\begin{eqnarray*}
\varepsilon^{3} f_k(R) 
& = & \varepsilon^5 R_{k}
\\&&+ i k \sum_{k' \in I_k \cap \vert k^\prime \vert \in \{1,3\}}(\varepsilon^3 R_{k-k'} + \varepsilon^2 Y_{k-k'})(\varepsilon Y_{k'}  + \varepsilon^2 R_{k'}) \\&&+ i k \sum_{\vert k- k' \vert = 1,3\; \cap \; \vert k^\prime \vert \notin \{1,3\}} (\varepsilon Y_{k-k'}  + \varepsilon^2 R_{k-k'})(\varepsilon^3 R_{k'} + \varepsilon^2 Y_{k'})\\&&+i k \sum_{\vert k- k' \vert = 1,3\; \cap \; \vert k^\prime \vert \in \{1,3\}} \big( \varepsilon^2 R_{k-k'}(\varepsilon^2 R_{k'}+\varepsilon Y_{k'}) +\varepsilon Y_{k-k'} \varepsilon^2 R_{k'} \big)
\\ && + i k \sum_{k' \in I_k \cap \vert k^\prime \vert \notin \{1,3\}}(\varepsilon^3 R_{k-k'} + \varepsilon^2 Y_{k-k'})(\varepsilon^3 R_{k'} + \varepsilon^2 Y_{k'})+   \textrm{Res}_r(k,t),
%
\end{eqnarray*}
and the reduced residual 
\begin{eqnarray*}
\textrm{Res}_r(\pm 2,t)  & = & - \varepsilon^4 \partial_T A_{\pm 2}+\varepsilon^4 Y_{\pm 2},\\
\textrm{Res}_r(\pm 4,t)  & = & - \varepsilon^4 \partial_T A_{\pm 4}+\varepsilon^4 Y_{\pm 4},\\
\textrm{Res}_r(\pm 6,t)  & = & - \varepsilon^4 \partial_T A_{\pm 6}+\varepsilon^4 Y_{\pm 6},
\end{eqnarray*} since $A_{\pm 2}, A_{\pm 4}$ and $A_{\pm 6}$ satisfy \eqref{A2eq}, \eqref{A4eq} and \eqref{A6eq} respectively. For all other $ k \in \Z $, it holds
$ 
\textrm{Res}_r(k,t) = 0 $.

\section{Estimates for the reduced residual}
\label{secresidual}

In this section we estimate the reduced residual
terms $
\textrm{Res}_r(\pm 2,t) $, $
\textrm{Res}_r(\pm 4,t) $, and $
\textrm{Res}_r(\pm 6,t) $. As a prototype example we consider 
$ \textrm{Res}_r(6,t) $.
%
From \eqref{A2eq}
we obtain  
\begin{eqnarray*}
\partial_T A_6 & = & - \frac{6i}{\lambda(6)} \partial_T((\widehat{Z}_3+A_3)^2) = 
- \frac{12i}{\lambda(6)} (\widehat{Z}_3+A_3)(\partial_T \widehat{Z}_3+\partial_T A_3) 
 \\ 
 & = & 
 -\frac{12i}{\lambda(6)} (\widehat{Z}_3+A_3)(\varepsilon^{-2}\partial_t \widehat{Z}_3+\partial_T A_3) 
\\ &  = &
- \frac{12i}{\lambda(6)}(\widehat{Z}_3+A_3)( \varepsilon^{-2} \alpha_3  \xi_3) \\ 
 && - \frac{12i}{\lambda(6)} (\widehat{Z}_3+A_3)
 (A_3 +  \widehat{Z}_3  
+ 6i (  \widehat{Z}_2  - \frac{1}{\lambda(2)} (2i(  \widehat{Z}_1  + A_1)^2 \\ && \nonumber \qquad  + 4i(  \widehat{Z}_3  + A_3)( \widehat{Z}_{-1}  + A_{-1}) ))
 ( \widehat{Z}_{1}  + A_{1})
\\ &&  \nonumber \qquad+ 6i (  \widehat{Z}_6  - \frac{6 i}{\lambda(6)}(  \widehat{Z}_3  + A_3)^2)
 ( \widehat{Z}_{-3}  + A_{-3}) 
+ 6 i \widehat{Z}_0 (\widehat{Z}_3+A_3)
\\ && \nonumber \qquad + 6i (  \widehat{Z}_4  - \frac{8 i}{\lambda(4)}  (  \widehat{Z}_3  + A_3)( \widehat{Z}_{1}  + A_{1}))
 ( \widehat{Z}_{-1}  + A_{-1})).
\end{eqnarray*}
Hence, most terms are classical functions but the term 
$$ 
 -\frac{12i}{\lambda(6)}(\widehat{Z}_3+A_3)( \varepsilon^{-2} \alpha_3  \xi_3) 
$$ 
needs some extra treatment. Therefore, we write 
$$
\textrm{Res}_r(6,t)  = \textrm{Res}_r(6,a,t) + \textrm{Res}_r(6,b,t),
$$
where 
$$ 
\textrm{Res}_r(6,b,t) =+ \varepsilon^{4}  \frac{12i}{\lambda(6)}(\widehat{Z}_3+A_3)( \varepsilon^{-2} \alpha_3  \xi_3) .
$$ 
In Section \ref{secerrorestimates} these residual terms appear as 
$$
\int_0^t e^{\lambda(6) (t-\tau)}
\varepsilon^{-3} \textrm{Res}_r(6,\tau) d\tau .
$$  
Under the assumption that 
\begin{equation} \label{cond3}
\sup_{T \in [0,T_0]} ( |A_1(T)| + |A_3(T)| + \sum_{j \in \{ 0,1,2,3,4,6\} }|Z_{j}(T/\varepsilon^2)|) \leq C_{\psi} = \mathcal{O}(1)
\end{equation}
we obviously have 
\begin{equation} \label{resesti2a}
\sup_{t \in [0,T_0/\varepsilon^2]}|\int_0^t e^{\lambda(6) (t-\tau)}
\varepsilon^{-3} \textrm{Res}_r(6,a,\tau) d\tau | \leq C \varepsilon.
\end{equation}
Less obvious is the estimate for 
$$ 
s_1(t) = |\int_0^t e^{\lambda(6) (t-\tau)}
  \varepsilon^{-3} \textrm{Res}_r(6,b,\tau) d\tau| .
$$ 
Recalling that $ \alpha_k = \mathcal{O}(\varepsilon^2)$ for $ k = \pm 1,\pm 3 $, we find 
\begin{eqnarray*}
s_1(t) & = &
   | \int_0^t e^{\lambda(6) (t-\tau)} \varepsilon^{-1} \frac{12i}{\lambda(6)}(\widehat{Z}_3+A_3) \alpha_3  \xi_3(\tau)d\tau|\\
   & = & 
   | \int_0^t e^{\lambda(6) (t-\tau)} \varepsilon  \frac{12i}{\lambda(6)}(\widehat{Z}_3+A_3)  \xi_3(\tau)d\tau|
\\& = & | \int_0^t e^{\lambda(6) (t-\tau)} \varepsilon  \frac{12i}{\lambda(6)}(\widehat{Z}_3+A_3) (\tau)  d\widehat{W}(3,\tau)|.
\end{eqnarray*}
Using \cite[(5.1.4)]{Arnold} again, we find 
\begin{eqnarray*}
\lefteqn{P(\sup_{\tau \in [0,t]}|s_1(\tau)| \geq c_k)}
\\& \leq & c_k^{-2} \int_0^t e^{2 \lambda(6) (t-\tau)} 
( \frac{12 \varepsilon}{\lambda(6)}  )^2  (E |Z_3(\tau)|^2+E |A_3(\tau)|^2) d\tau
\\& \leq &  \frac{144}{\lambda(6)^2} \varepsilon^2  c_k^{-2}  \int_0^t e^{2\lambda(6) (t-\tau)} 
  d\tau 
(\sup_{\tau \in [0,t]} (E |Z_3(\tau)|^2)+\sup_{\tau \in [0,t]} (E |A_3(\tau)|^2))
\\ & \leq &  \frac{144}{\lambda(6)^2}  \frac{\varepsilon^2}{2 \lambda(6) c_k^2} 
(E \sup_{\tau \in [0,t]} ( |Z_3(\tau)|^2)+E \sup_{\tau \in [0,t]} ( |A_3(\tau)|^2)).
\end{eqnarray*}
Using \cite[(5.1.5)]{Arnold} gives the estimate for $ E \sup_{\tau \in [0,t]} ( |Z_1(\tau)|^2) $. Similarly we can  estimate $ E \sup_{\tau \in [0,t]} ( |A_1(\tau)|^2) $.

Almost line for line we find the same estimates for $ j = 4,6 $.
Summarizing the estimates yields that for all $ \delta> 0 $ there 
exists a $ C > 0 $ such for all 
$ \varepsilon \in (0,1) $ we have that 
\begin{equation} \label{resesti2a}
P(\sup_{t \in [0,T_0/\varepsilon^2]}|\int_0^t e^{\lambda(j) (t-\tau)}
\varepsilon^{-3} \textrm{Res}_r(j,\tau) d\tau | \leq C_{res} \varepsilon) > 1-\delta
\end{equation}
for $ j = 2,4,6 $.

\section{The  error estimates}
\label{secerrorestimates}

In order to estimate the error made by the stochastic Landau approximation
we use the variation of constant of formula and write the
equations \eqref{r1eq}, \eqref{r3eq} and \eqref{rkeq} for the error as
\begin{eqnarray*}
R_1(t) & = & e^{\lambda(1) t} R_1(0) + \int_0^t e^{\lambda(1) (t-\tau)} f_1(\tau) d \tau, \\
R_3(t) & = & e^{\lambda(3) t} R_3(0) + \int_0^t e^{\lambda(3) (t-\tau)} f_3(\tau) d \tau, \\ 
R_k(t) & = & e^{\lambda(k) t} R_k(0) + \int_0^t e^{\lambda(k) (t-\tau)} 
f_k(\tau)
 d \tau. 
\end{eqnarray*}
We distinguish the error function for critical and stable modes by setting
$$
R_c = (R_{-3},R_{-1},R_1,R_3) , \qquad R_s= (\ldots, R_{-4},R_{-2},R_0,R_2,R_4,\ldots)
$$
and define the quantities
$$
S_c(t) = \sup_{0 \leq \tau \leq t} \| R_c(\tau) \|_{\ell^2_r}
\qquad \textrm{and} \qquad
S_s(t) = \sup_{0 \leq \tau \leq t} \| R_s(\tau) \|_{\ell^2_r},
$$
with $ r \in (1/2,3) $.
The following estimates hold under the assumption that 
\begin{equation} \label{cond1}
\sum_{k \in \{ 1,2,3,4,6\} }\sup_{T \in [0,T_0]}  |A_k(T)|  \leq C_{\psi} 
\end{equation}
and
\begin{equation} \label{cond2}
\sup_{t \in [0,T_0/\varepsilon^2]} \sum_{k\in \Z}  |\widehat{Z}_k(t) |^2 (1+|k|^2)^r \leq C_{Z}^2.
\end{equation}

i) Using \eqref{convul} and $ \lambda(1) = 0 $ we estimate 
\begin{eqnarray*}
&& |\int_0^t \varepsilon^2 R_1(\tau) d\tau | \leq C \int_0^t \varepsilon^2 S_c(\tau) d\tau, \\ 
&&  |\int_0^t 
 i  (\sum_{1- k' = \pm 1,\pm 3 } \big(R_{1-k'}(\varepsilon^3 R_{k'} +\varepsilon^2 Y_{k'}) +\varepsilon Y_{1-k'} \varepsilon R_{k'} \big))
 (\tau) d \tau | \\ && \qquad
\leq C \int_0^t  \varepsilon^3 S_c(\tau) S_s(\tau) + \varepsilon^2 S_c(\tau) +\varepsilon^2 S_s(\tau)d\tau, \\
&&  |\int_0^t   i  \sum_{k' \in I_1 \; \cap \;\vert k^\prime \vert \in \{1,3\}} \big(\varepsilon R_{1-k'} (\varepsilon^2 R_{k'}+\varepsilon Y_{k'}) +Y_{1-k'} \varepsilon^2 R_{k'} \big)(\tau) d \tau | \\&& \qquad \leq C \int_0^t  \varepsilon^3 S_c(\tau) S_s(\tau) + \varepsilon^2 S_s(\tau) +\varepsilon^2 S_c(\tau) d\tau, 
\\ 
&& |\int_0^t 
 i \sum_{k' \in I_1 \; \cap \;\vert k^\prime \vert \notin \{1,3\}}(\varepsilon^3 R_{1-k'} + \varepsilon^2 Y_{1-k'})(\varepsilon R_{k'} +  Y_{k'})) (\tau) d \tau | \\ &&  \qquad
\leq C \int_0^t \varepsilon^2  + \varepsilon^3  S_s(\tau) + \varepsilon^4 S_s^2(\tau) d\tau, \\
\\ 
&& |\int_0^t
  \varepsilon^{-2} \textrm{Res}_r(1,t)  (\tau) d \tau | = 0.
\end{eqnarray*}
ii) In exactly the same way we estimate the terms in the equations for $ R_3 $.

\noindent
iii) The estimates for  the stable part are more advanced. 
They are fundamentally based on the fact  that  $ \lambda(k) < 0 $ for $ k \neq \pm 1,\pm 3 $. 
We start with 
\begin{eqnarray*}
(\sum_{k \in \Z}|\int_0^t e^{\lambda(k) (t-\tau)} \varepsilon^2 R_k(\tau) d\tau|^2 (1+k^2)^r)^{1/2}
\leq C (\varepsilon^2 S_s(t)) .
\end{eqnarray*}
Next we find
\begin{eqnarray*}
	&&(| \int_0^t e^{\lambda(k) (t-\tau)}
	i k ( \sum_{k' \in I_k \cap \vert k^\prime \vert \in \{1,3\}}(\varepsilon R_{k-k'} + Y_{k-k'})(Y_{k'}  + \varepsilon R_{k'})(\tau)d\tau|^2 (1+k^2)^r)^{1/2} \\ && \leq  C ( \varepsilon S_s(t)+\varepsilon^2 S_s(t) S_c(t)+\varepsilon S_c(t)+ 1)
\end{eqnarray*}
%
%
%
and
\begin{eqnarray*}
	&&(| \int_0^t e^{\lambda(k) (t-\tau)}
	i k( \sum_{\vert k- k' \vert = 1,3\; \cap \; \vert k^\prime \vert \notin \{1,3\}}  (Y_{k-k'}  + \varepsilon R_{k-k'})(\varepsilon R_{k'} +Y_{k'}) )
	(\tau) d\tau|^2 (1+k^2)^r)^{1/2} 
	\\ && \leq  C (\varepsilon S_s(t)+1+\varepsilon^2 S_c(t)S_s(t)+\varepsilon S_c(t)).
\end{eqnarray*}
%
%
%
Moreover, we have 
\begin{eqnarray*}
	&&(| \int_0^t e^{\lambda(k) (t-\tau)}
	i k (\sum_{\vert k- k' \vert = 1,3\; \cap \; \vert k^\prime \vert \in \{1,3\}} \big( R_{k-k'}(\varepsilon R_{k'}+Y_{k'}) + Y_{k-k'} R_{k'} \big)
	)
	(\tau) d\tau|^2 (1+k^2)^r)^{1/2} \\&& \qquad
	\leq  C (\varepsilon S_c^2(t)+2S_c(t)).
\end{eqnarray*}
%
%
%
%
%
Furthermore, we estimate
\begin{eqnarray*}
	&& (| \int_0^t e^{\lambda(k) (t-\tau)}
	 i k (\sum_{k' \in I_k \cap \vert k^\prime \vert \notin \{1,3\}}(\varepsilon R_{k-k'} + Y_{k-k'})(\varepsilon^2 R_{k'} + \varepsilon Y_{k'})) (\tau) d\tau|^2 (1+k^2)^r)^{1/2} 
	\\ && \leq  C \varepsilon(\varepsilon^2 S_s^2(t) +\varepsilon S_s(t) +1).
\end{eqnarray*}
%
%
%
%

\medspace

Recalling that $ \textrm{Res}_r(k,t) = 0  $ for $ k \neq \pm 2 ,\pm 4,\pm 6$,  \eqref{resesti2a} yields
\begin{eqnarray*}
&& (| \int_0^t e^{\lambda(k) (t-\tau)}
 \varepsilon^{-3} \textrm{Res}_r(k,\tau) d\tau|^2 (1+k^2)^r)^{1/2} 
 \leq  C \varepsilon .
\end{eqnarray*}
Now we have all ingredients to establish the  bound for the error.
For all $  t \in [0, T_0/ \varepsilon^2] $ we have 
\begin{eqnarray*}
S_c(t)  & \leq  &  S_c(0) + 
 \int_0^t  C \varepsilon^2 (1+3S_c(\tau) + 2S_s(\tau)+ 2\varepsilon S_c(\tau) S_s(\tau)+ \varepsilon S_s(\tau)  +\varepsilon^2 S_s(\tau)^2 )d\tau
\\S_s(t)  & \leq  &  S_s(0) + 2C(S_c(t)+1)+ C\varepsilon[(M_c+\varepsilon M_s)^2+2M_c(t)+2M_s(t) +2\varepsilon M_s(t) +1],
\end{eqnarray*}
%

From the second inequality we then obtain:
There exists an $ \varepsilon_1 > 0 $ such that 
for all $ \varepsilon \in (0,\varepsilon_1) $ we have 
\begin{equation}\label{gries}
S_s(t) \leq  (S_s(0)+1) + 2C(S_c(t)+1),
\end{equation}
as long as $ S_c(t) \leq M_c $ and $ S_s(t) \leq M_s $ with $ M_c $ and $ M_s $ defined below,
if $ 0 < \varepsilon_1 \ll 1$  is chosen so small that 

\begin{equation} \label{zweibr1}
	C\varepsilon[(M_c+\varepsilon M_s)^2+2M_c(t)+2M_s(t) +2\varepsilon M_s(t) +1] \leq 1
\end{equation}
for all $ \varepsilon \in (0,\varepsilon_1) $.
Inserting \eqref{gries} into the inequality for $ S_c(t) $ yields 
$$ 
S_c(t)  \leq   \beta_0 + \varepsilon^2 \int_0^t \beta_1 S_c(\tau)   d\tau,
$$
where
\begin{eqnarray*}
\beta_0 &= & S_c(0)+1
+ C T_0(1+2S_s(0)+2+4C),  \\
\beta_1&= & 
  C(3+4C)  ,
\end{eqnarray*}
if we choose 
$ 0 < \varepsilon_2 \ll 1 $ so small that for given $ M_c $ and $ M_s $ we have 
\begin{equation}
\label{saarbr1}
 CT_0 (\varepsilon S_s(0)+\varepsilon+\varepsilon M_c M_s+ \varepsilon^2 M_s^2)  \leq 1 
\end{equation}
for all $ \varepsilon \in (0,\varepsilon_2) $.
Using Gronwall's inequality immediately 
yields 
$$ 
S_c(t) \leq  \beta_0 e^{\beta_1 T_0}=:M_c, 
$$ 
and due to \eqref{gries} we define  
$$
M_s:= (S_s(0)+1) + 2C(S_c(t)+1).
$$
For these $ M_c $ and $ M_s $ we define 
 $ \varepsilon_0 = \min(\varepsilon_1,\varepsilon_2) $ where 
$ \varepsilon_1 > 0 $ is chosen so small that for all $ \varepsilon \in (0,\varepsilon_1) $ condition 
\eqref{zweibr1} is satisfied
and 
$ \varepsilon_2 > 0 $  so small that for all $ \varepsilon \in (0,\varepsilon_2) $ condition 
\eqref{saarbr1} is satisfied.

Since the above finitely many 
conditions \eqref{cond3}, \eqref{resesti2a}, \eqref{cond1}, and \eqref{cond2} can be satisfied with a full probability independently of the small parameter $ 0 < \varepsilon \ll 1 $, cf. \eqref{Zesti},  \eqref{rem42}, and \eqref{resesti2a}, if the $ \mathcal{O}(1) $-bounds $ C_{res}$, $ C_{\psi} $, and $ C_Z^2 $ go to infinity, 
we have: 

For all $ \delta > 0 $  there exist $ \varepsilon_0 > 0 $, $ C_2 > 0 $  such that for all $ \varepsilon \in (0,\varepsilon_0) $ there are solutions 
$ (R_k)_{k \in \Z} \in C([0,T_0/\varepsilon^2],\ell^2_r) $
of  \eqref{r1eq} and \eqref{rkeq} with 
\begin{eqnarray*}
&& \mathcal{P}(\sup_{t \in [0,T_0/\varepsilon^2]} 
(\| (R_k(t))_{k \in \Z} \|_{\ell^2_r}   
\leq C_2 \varepsilon^2)   >   1 - \delta.
\end{eqnarray*}
Sobolev's embedding theorem yields the statement of Theorem \ref{mainth}. \qed



\section{Discussion}
\label{secdisc}

In the previous chapters, we brought the derivation and justification of stochastic amplitude equations for pattern-forming systems with additive white noise  closer to the deterministic situation. 
It will be straightforward to transfer the presented approach to more complicated pattern forming systems with periodic boundary conditions.

For notational simplicity we restricted the previous analysis to the duKS equation 
\eqref{org1}. Another simple example is a  system of coupled 
Swift-Hohenberg-Kuramoto-Shivashinsky (SHKS) equations 
\begin{eqnarray*}
\partial_t u  &  = &  -(k_1^2+\partial_x^2)^2 u + \alpha u + f_u(u,v,\partial_x u ,\partial_x v)+ \xi_u, \\
\partial_t v  &  = &  -(k_2^2+\partial_x^2)^2 v + \alpha v + f_v(u,v,\partial_x u ,\partial_x v) + \xi_v
\end{eqnarray*}
where $x, \, u=u(x,t), v=v(x,t) \in \mathbb{R}, \, t \geq 0 $, bifurcation parameter $ \alpha \in  \mathbb{R} $, smooth nonlinearities $ f_u,f_v: \R^4 \to \R $
with $ f_{u,v}(u,v,\partial_x u ,\partial_x v) = \mathcal{O}(|u|^2+ |v|^2+|\partial_x u|^2 + |\partial_x v|^2)$ for $ u,v,\partial_x u,\partial_x v \to 0 $,
and noise $ \xi_u= \xi_u(x,t) $ and $ \xi_v= \xi_v(x,t) $.
As before we have to satisfy the condition \eqref{Zesti1a}
which holds if $ 2(s+r) < 3 $
or equivalently
$0 \leq  s+ r < 3/2 $
since the eigenvalues satisfy $  \lambda_{1,2}(k) = -(k_j^2 - k^2)^2   \sim -k^4  $ for $|k| \to \infty $.
Hence, the subsequent error equations can be solved in the space $ \ell^2_r $ with 
$ 1/2 < r < 1 $.

If this instability mechanism is considered for reaction diffusion systems 
$$
\partial_t U = D \partial_x^2 U + f(U)
$$
with diffusion matrix $ D \in \R^{d \times d} $, $ U(x,t) \in \R^d $ and $ f : \R^d \to \R^d $ 
a smooth mapping, the condition   \eqref{Zesti1a} can only be satisfied 
if the white noise, i.e., $ \alpha_k = \mathcal{O}(1) $ for $ |k | \to \infty $
is replaced by colored noise, i.e., $ \alpha_k = \mathcal{O}(|k|^{-\theta}) $ for $ |k | \to \infty $ and $ \theta > 0 $. Including $ \alpha_k $ in \eqref{Zesti1a} gives 
the condition 
\begin{equation} \label{Zesti1ab}
\sum_{k\in \Z\setminus I_u} \frac{C_s^{2} \alpha_k^2 (1+k^2)^{(s+r)} }{2 \lambda(k)} = \mathcal{O}(1)< \infty 
\end{equation} 
which is the case if $ 2(s+r-\theta) < 1 $
or equivalently
$0 \leq  s+ r < 1/2 + \theta  $
since $  \lambda(k) \sim -k^2 $ for $|k| \to \infty $.
Hence, the subsequent error equations can be solved in the space $ \ell^2_r $ with 
$ 1/2 < r < \theta $.

The same analysis for B\'enard's problem with two layered fluids is more
involved and 
will be the subject of future research.

Another future goal is to use  the presented approach  to transfer
more complicated deterministic bifurcation scenarios to the stochastic case.
As a first step in this direction  the 
$ 2 \pi $-periodicity   will be given up 
and the original system will be considered on a domain of size 
$ \mathcal{O}(1/\varepsilon) $ with periodic boundary conditions.
Similar to \cite{BHP05} the stochastic Landau equation will be replaced by 
a stochastic Ginzburg-Landau equation.

Last but not least, it is worth mentioning that the stable part in \eqref{futureWork} can be chosen of order $\varepsilon$ instead of $\varepsilon^2$, as pointed out by Dirk Blömker. Letting $c_k=\varepsilon$ for $k\neq \pm 1,3$, would deliver an improvement to the present results. To provide however a more rigorous analysis is part of our future work.

\bibliographystyle{plain}
\bibliography{GLbib}

\end{document}